\documentclass[12pt]{article}
\usepackage{amssymb,amsmath}
\usepackage{amsthm}
\usepackage{cases}
\usepackage{amsfonts}
\usepackage{graphicx}
\usepackage{makecell}
\usepackage{cite,color,xcolor}
\usepackage[left=1.90cm,right=1.90cm,top=1.80cm,bottom=1.80cm]{geometry}
\usepackage[colorlinks,citecolor=blue,urlcolor=blue]{hyperref}
\usepackage[utf8]{inputenc}

\newtheorem{theorem}{Theorem}[section]
\newtheorem{corollary}{Corollary}[section]
\newtheorem{lemma}{Lemma}[section]
\newtheorem{proposition}{Proposition}[section]

\newtheorem{definition}{Definition}[section]
\newtheorem{remark}{Remark}[section]

\usepackage{txfonts}
\newcommand{\bal}{\begin{align}}
\newcommand{\bbal}{\begin{align*}}
\newcommand{\beq}{\begin{equation}}
\newcommand{\eeq}{\end{equation}}
\newcommand{\bca}{\begin{cases}}
\newcommand{\eca}{\end{cases}}

\newcommand{\pa}{\partial}
\newcommand{\fr}{\frac}

\newcommand{\De}{\Delta}

\newcommand{\ep}{\varepsilon}
\newcommand{\dd}{\mathrm{d}}

\newcommand{\R}{\mathbb{R}}

\newcommand{\les}{\lesssim}

\newcommand{\f}{\left}
\newcommand{\g}{\right}

\begin{document}
\title{Nonexistence of H\"{o}lder continuous solution for the Camassa-Holm equation in Besov spaces }
\author{Yanghai Yu$^{1,}$\footnote{E-mail: yuyanghai214@sina.com(Corresponding author); lijinlu@gnnu.edu.cn; mathzwp2010@163.com}, Jinlu Li$^{2}$  and Weipeng Zhu$^{3}$\\
\small $^1$ School of Mathematics and Statistics, Anhui Normal University, Wuhu 241002, China\\
\small $^2$ School of Mathematics and Computer Sciences, Gannan Normal University, Ganzhou 341000, China\\
\small $^3$ School of Mathematics, Foshan University, Foshan, Guangdong 528000, China}

\date{}

\maketitle\noindent{\hrulefill}

{\bf Abstract:}
In the paper, we show that the continuity of the solution can not be improved to the H\"{o}lder continuity. Precisely speaking, the solution of the Camassa-Holm equation belongs to $\mathcal{C}([0,T];B^s_{p,r})$ but not to $\mathcal{C}^\alpha([0,T];B^s_{p,r})$ with any $\alpha\in(0,1)$. To the best of our knowledge, our work is the first one addressing the issue on the failure of H\"{o}lder continuous in time of solution to the classical Camassa-Holm equation.
As a by-product, we establish the ill-posedness for the Camassa-Holm equation in $B^s_{p,\infty}(\mathbb{R})$ with $s>\max\big\{1+1/p, 3/2\big\}$ with $p\in[1,\infty]$ by proving the solution map to the Camassa-Holm equation starting from $u_0$ is discontinuous at $t = 0$ in $B^s_{p,\infty}(\mathbb{R})$.


\vskip0mm\noindent{\hrulefill}

\section{Introduction}

In this paper, we are concerned with the Cauchy problem for the classical Camassa-Holm (CH) equation
\begin{equation}\label{0}
\begin{cases}
u_t-u_{xxt}+3uu_x=2u_xu_{xx}+uu_{xxx}, \; &(x,t)\in \R\times\R^+,\\
u(x,t=0)=u_0(x),\; &x\in \R.
\end{cases}
\end{equation}
Here the scalar function $u = u(t, x)$ stands for the fluid velocity at time $t\geq0$ in the $x$ direction.
Setting $\Lambda^{-2}=(1-\pa^2_{xx})^{-1}$, then $\Lambda^{-2}f=G*f$ where $G(x)=\fr12e^{-|x|}$ is the kernel of the operator $\Lambda^{-2}$. Thus, we can transform the CH equation \eqref{0} equivalently into the following transport type equation
\begin{equation}\label{CH}
\begin{cases}
\partial_tu+u\pa_xu=\mathbf{P}(u):=-\pa_x\Lambda^{-2}\left(u^2+\fr12(\pa_xu)^2\right), \; &(x,t)\in \R\times\R^+,\\
u(x,t=0)=u_0(x),\; &x\in \R.
\end{cases}
\end{equation}
The CH equation \eqref{0} was first proposed in the context of hereditary symmetries studied in \cite{Fokas} and then was derived explicitly as a water wave equation by Camassa-Holm \cite{Camassa}. Many aspects of the mathematical beauty of the CH equation have been exposed over the last two decades. Particularly,  \eqref{0} is completely integrable \cite{Camassa,Constantin-P} with a bi-Hamiltonian structure \cite{Constantin-E,Fokas} and infinitely many conservation laws \cite{Camassa,Fokas}. Also, it admits exact peaked soliton solutions (peakons) of the form $u(x,t)=ce^{-|x-ct|}\;(c>0)$, which are orbitally stable \cite{Constantin.Strauss}. Another remarkable feature of the CH equation is the wave breaking phenomena: the solution remains bounded while its slope becomes unbounded in finite time \cite{Constantin}.
After the CH equation was derived physically in the context of water waves, there is a large amount of literature devoted to studying the well-posedness of the Cauchy problem \eqref{0} in Besov spaces. Danchin \cite{d1,d3} proved the local existence and uniqueness of strong solutions to \eqref{0} with initial data in $B^s_{p,r}$ if $(p,r)\in[1,\infty]\times[1,\infty), s>\max\big\{1+1/p, 3/2\big\}$ and $B^{3/2}_{2,1}$.  Meanwhile, he \cite{d1} only obtained the continuity of the solution map of \eqref{0} with respect to initial data in the space $\mathcal{C}([0, T ];B^{s'}_{p,r})$ with any $s'<s$. Li-Yin \cite{Li-Yin1} proved the continuity of the solution map of \eqref{0} with respect to initial data in the space $\mathcal{C}([0, T];B^{s}_{p,r})$ with $r<\infty$. In particular, they \cite{Li-Yin1} proved that the solution map of \eqref{0} is weak continuous with respect to initial data $u_0\in B^s_{p,\infty}$. For the endpoints, Danchin \cite{d3} obtained that the data-to-solution map is not continuous by using peakon solution, which implies the ill-posedness of \eqref{0} in $B^{3/2}_{2,\infty}$.
Guo-Liu-Molinet-Yin \cite{Guo-Yin} showed the ill-posedness of \eqref{0} in $B_{p,r}^{1+1/p}(\mathbb{R}\;\text{or}\; \mathbb{T})$ with $(p,r)\in[1,\infty]\times(1,\infty]$ by proving the norm inflation. Guo-Ye-Yin \cite{Guo} and Li-Yu-Zhu \cite{24jmfm,26ampa} obtained the ill-posedness for the CH equation and the Novikov equation in $B^{1}_{\infty,1}$ by proving the norm inflation, respectively. Ye-Yin-Guo \cite{Ye} obtained the local well-posedness for the Camassa-Holm type equation in $B^s_{p,r}$ with $s>1+\frac1p,  (p,r)\in[1,+\infty]\times[1,+\infty)$ or $s=1+\frac1p, \ (p,r)\in[1,+\infty)\times \{1\}$.
Assume that $u_0\in B^s_{p,r}$ with $s>1+\frac1p,  (p,r)\in[1,+\infty]\times[1,+\infty)$ or $s=1+\frac1p, \ (p,r)\in[1,+\infty)\times \{1\}$, it is known that there exists a solution $u\in \mathcal{C}([0,T];B^s_{p,r})$ for the Camassa-Holm equation. Naturally, we may wonder whether or not the solution $u$ can belong to $\mathcal{C}^\alpha([0,T];B^s_{p,r})$ with some $\alpha\in(0,1)$. Therefore, we are interested in the following
{\bf Question}:
$$u_0\in B^s_{p,r}\;\Rightarrow\; \exists 1\;u\in \mathcal{C}([0,T];B^s_{p,r})\;\overset{?}{\Rightarrow}\; u\in\mathcal{C}^\alpha([0,T];B^s_{p,r}) \quad\text{with}\; \alpha\in(0,1).$$

Assume that  the initial data $u_0$ has more regularity such that $u_0\in B^{s'}_{p,r}$ for some $s'>s$, by the interpolation argument, we can deduce that $u\in \mathcal{C}^\alpha([0,T];B^s_{p,r})$ with $\alpha=s'-s$.
In this paper, we will show that that there exists
initial data $u_0\in B^s_{p,r}$ such that the corresponding solution of the Camassa-Holm equation can not belong to $\mathcal{C}^\alpha([0,T];B^s_{p,r})$ with any $\alpha\in(0,1)$. Now let us state our main result of this paper.
\begin{theorem}\label{th1}
Assume that $(s,p,r)$ satisfies that
\bal\label{cond}
s>1+\frac1p, \ (p,r)\in[1,\infty]\times[1,\infty) \quad\text{or}\quad s=1+\frac1p, \ (p,r)\in[1,\infty)\times \{1\}.
\end{align}
For any $\alpha\in(0,1)$, there exists $u_0\in B^s_{p,r}(\R)$  such that the data-to-solution map $u_0\mapsto \mathbf{S}_{t}(u_0)\in \mathcal{C}([0,T];B^s_{p,r})$ of the Cauchy problem \eqref{CH}
satisfies
\bbal
\limsup_{t\to0^+}\frac{\|\mathbf{S}_{t}(u_0)-u_0\|_{B^s_{p,r}}}{t^\alpha}=+\infty.
\end{align*}
\end{theorem}
In our recent papers \cite{jee}, we established the ill-posedness for the Camassa-Holm equation in $B^s_{p,\infty}(\mathbb{R})$ with $s>\max\big\{1+1/p, 3/2\big\}$ with $p\in[1,\infty]$ by proving the solution map to the CH equation starting from $u_0$ is discontinuous at $t = 0$ in the metric of $B^s_{p,\infty}(\mathbb{R})$. As a by-product of Theorem \ref{th1}, we have
\begin{corollary}\label{co1}
Assume that $(s,p,r)$ satisfies that
$s>1+\frac1p$ with $p\in[1,\infty]$.
There exits $u_0\in B^s_{p,\infty}(\R)$ and a positive constant $\ep_0$ such that the data-to-solution map $u_0\mapsto \mathbf{S}_{t}(u_0)$ of the Cauchy problem \eqref{CH}
satisfies
\bbal
\limsup_{t\to0^+}\|\mathbf{S}_{t}(u_0)-u_0\|_{B^s_{p,\infty}}\geq\ep_0.
\end{align*}
\end{corollary}
\begin{remark} Compared with the previous result in \cite{jee}, Corollary \ref{co1} is new since the regularity have been enlarged. In fact, based on some new observations and subtle estimations, we can deal with the case when $1+\frac{1}{p}<s\leq\frac{3}{2}$.
\end{remark}

\begin{remark}
We summary the local well-posedness/ill-posdeness results of the Camassa-Holm equation in Besov spaces. This can be seen clearly from the Table 1.
\end{remark}
\begin{table}[http]
      \centering
      \begin{tabular}{l|c|c}\hline
       References&Range&Results\\\hline
        \cite{d1,d3,Li-Yin1}&\makecell[c]{$s>\max\left\{1+1/p, 3/2\right\},(p,r)\in[1,\infty]\times[1,\infty)$ \\$s=3/2,p=2,r=1$}&LWP \\\hline
        \cite{Guo-Yin}&$s=1+1/p,(p,r)\in[1,\infty]\times(1,\infty]$&Norm inflation
        \\\hline
        \cite{Guo}&$s=1,p=\infty,r=1$& Norm inflation
          \\\hline
        \cite{Ye}&\makecell[c]{$s>1+1/p, \ (p,r)\in[1,\infty]\times[1,\infty)$\\ $s=1+1/p, \ (p,r)\in[1,\infty)\times \{1\}$ }&LWP
        \\\hline
        Theorem \ref{th1}&\makecell[c]{$s>1+1/p, \ (p,r)\in[1,\infty]\times[1,\infty)$\\ $s=1+1/p, \ (p,r)\in[1,\infty)\times \{1\}$ }& No H\"{o}lder continuos \\\hline
        Corollary \ref{co1}&$s>1+1/p, \ p\in[1,\infty],r=\infty$& Discontinuos \\\hline
        \end{tabular}
        \caption{Well/Ill-posedness of \eqref{CH} in $B^s_{p,r}$}
        \end{table}

\begin{remark}
Theorem \ref{th1} is true for the Torus case and also holds for the b-family equation
\begin{align*}
\begin{cases}
\pa_tu+uu_x=-\pa_x(1-\pa^2_x)^{-1}\left(\frac{b}{2}u^2
+\frac{3-b}{2}u^2_x\right), &\quad (t,x)\in \R^+\times\R,\\
u(0,x)=u_0(x), &\quad x\in \R,
\end{cases}
\end{align*}
and the Novikov equation
\begin{align*}
\begin{cases}
\pa_tv+v^2v_x=-(1-\pa^2_x)^{-1}\left(\frac12v_x^3
+\pa_x\left(\frac32vv^2_x+v^3\right)\right), \\
v(0,x)=v_0(x).
\end{cases}
\end{align*}
\end{remark}
\section{Preliminaries}\label{sec2}
{\bf Notation}\;
The notation $A\les B$ (resp., $A \gtrsim B$) means that there exists a harmless positive constant $c$ such that $A \leq cB$ (resp., $A \geq cB$).
Given a Banach space $X$, we denote its norm by $\|\cdot\|_{X}$. For $I\subset\R$, we denote by $\mathcal{C}(I;X)$ the set of continuous functions on $I$ with values in $X$. Sometimes we will denote $L^p(0,T;X)$ by $L_T^pX$.
Next, we will recall some facts about the Littlewood-Paley decomposition and the nonhomogeneous Besov spaces (see \cite{BCD} for more details).
Choose a radial, non-negative, smooth function $\vartheta:\R\mapsto [0,1]$ such that
 ${\rm{supp}} \,\vartheta\subset B(0, 4/3)$ and $\vartheta(\xi)\equiv1$ for $|\xi|\leq3/4$.
Setting $\varphi(\xi):=\vartheta(\xi/2)-\vartheta(\xi)$, then we deduce that $\varphi$ has the following properties
${\rm{supp}} \;\varphi\subset \left\{\xi\in \R: 3/4\leq|\xi|\leq8/3\right\}$ and $\varphi(\xi)\equiv 1$ for $4/3\leq |\xi|\leq 3/2$.

\begin{definition}[see \cite{BCD}]
For every $u\in \mathcal{S'}(\mathbb{R})$, the Littlewood-Paley dyadic blocks ${\Delta}_j$ are defined as follows
\begin{align*}
\Delta_ju=0,\; \text{if}\; j\leq-2;\quad
\Delta_{-1}u=\vartheta(D)u;\quad
\Delta_ju=\varphi(2^{-j}D)u,\; \; \text{if}\;j\geq0,
\end{align*}
The inhomogeneous low-frequency cut-off operator $S_j$ is defined by
$
S_ju=\sum_{q=-1}^{j-1}{\Delta}_qu.
$
\end{definition}
\begin{definition}[see \cite{BCD}]
Let $s\in\mathbb{R}$ and $(p,r)\in[1, \infty]^2$. The nonhomogeneous Besov space $B^{s}_{p,r}(\R)$ is defined by
\begin{align*}
B^{s}_{p,r}(\R):=\Big\{f\in \mathcal{S}'(\R):\;\|f\|_{B^{s}_{p,r}(\mathbb{R})}<\infty\Big\},
\end{align*}
where $\|f\|_{B^{s}_{p,\infty}(\mathbb{R})}=\sup_{j\geq-1}2^{sj}\|\Delta_jf\|_{L^p(\mathbb{R})}$ and for $1\leq r<\infty$,
$$\|f\|_{B^{s}_{p,r}(\mathbb{R})}=
\left(\sum_{j\geq-1}2^{sjr}\|\Delta_jf\|^r_{L^p(\mathbb{R})}\right)^{1/r}.$$
\end{definition}

We shall use Bony's decomposition \cite{BCD} in the nonhomogeneous context throughout this paper
\begin{align*}
uv={T}_{u}v+{T}_{v}u+{R}(u,v)\quad\text{with}\quad
{T}_{u}v=\sum_{j\geq-1}{S}_{j-1}u{\Delta}_jv \quad\mbox{and} \quad{R}(u,v)=\sum_{|j-k|\leq1}{\Delta}_ju{\Delta}_kv.
\end{align*}
\begin{lemma} [see \cite{BCD}] \label{leem1} Let $\left(s, t, s_{1}, s_{2}\right)\in\mathbb{R}^{4}$ and $\left(p, p_{1}, p_{2}, r, r_{1}, r_{2}\right)\in[1, \infty]^{6}$. Assume that
$$
\frac{1}{p} = \frac{1}{p_{1}}+\frac{1}{p_{2}} \leq 1, \quad \frac{1}{r} =\frac{1}{r_{1}}+\frac{1}{r_{2}} \leq 1, \quad s_{1}+s_{2}>0, \quad t<0 .
$$
There exists a constant $C$ such that
\bbal
&\left\|T_{u} v\right\|_{B_{p, r}^{s}\left(\mathbb{R}\right)} \leq C^{|s|+1}\|u\|_{L^{\infty}\left(\mathbb{R}\right)}\|v\|_{B_{p, r}^{s}\left(\mathbb{R}\right)},\\
&\left\|T_{u} v\right\|_{B_{p, r}^{s+t}\left(\mathbb{R}\right)} \leq \frac{C^{|s+t|+1}}{-t}\|u\|_{B_{\infty, r_{1}}^{t}\left(\mathbb{R}\right)}\|v\|_{B_{p, r_{2}}^{s}\left(\mathbb{R}\right)},\\
&\|R(u, v)\|_{B_{p, r}^{s_{1}+s_{2}}\left(\mathbb{R}\right)} \leq \frac{C^{\left|s_{1}+s_{2}\right|+1}}{s_{1}+s_{2}}\|u\|_{B_{p_{1}, r_{1}}^{s_{1}}\left(\mathbb{R}\right)}\|v\|_{B_{p_{2}, r_{2}}^{s_{2}}\left(\mathbb{R}\right)}.
\end{align*}
\end{lemma}
Finally, we give some important properties which will be also often used throughout the paper.
\begin{lemma}[see \cite{BCD}]\label{lm2}
Let $(p,r)\in[1, \infty]^2$ and $s>\max\big\{1+\frac1p,\frac32\big\}$. Then we have
\bbal
&\|uv\|_{B^{s-2}_{p,r}(\R)}\leq C\|u\|_{B^{s-2}_{p,r}(\R)}\|v\|_{B^{s-1}_{p,r}(\R)}.
\end{align*}
\end{lemma}

\begin{lemma}[see \cite{BCD}]\label{lm3}
For $(p,r)\in[1, \infty]^2$, $B^{s-1}_{p,r}(\R)$ with $s>1+\frac{1}{p}$ or $s=1+\frac{1}{p},\ r=1$ is an algebra. Moreover, for any $u,v \in B^{s-1}_{p,r}(\R)$ with $s>1+\frac{1}{p}$ or $s=1+\frac{1}{p},\ r=1$ , we have
\bbal
&\|uv\|_{B^{s-1}_{p,r}(\R)}\leq C\|u\|_{B^{s-1}_{p,r}(\R)}\|v\|_{B^{s-1}_{p,r}(\R)}.
\end{align*}
\end{lemma}
\begin{remark}\label{re5} Let $(p,r)\in[1, \infty]^2$ and $s>\max\big\{1+\frac1p,\frac32\big\}$, using Lemmas \ref{lm2}-\ref{lm3}, we have
\begin{align*}
\|\mathbf{P}(u)-\mathbf{P}(v)\|_{B_{p, r}^{s-1}} & \leq C\|u-v\|_{B_{p, r}^{s-1}}\f(\|u\|_{B_{p, r}^{s}}+\|v\|_{B_{p, r}^{s}}\g).
\end{align*}
\end{remark}

\begin{lemma}[see \cite{BCD}]\label{lm5}
Let $p\in[1,\infty]$ and $s>0$. There exists
a constant $C$, depending continuously on $p$ and $s$, such that
\bbal
\f\|2^{j s}\|[\Delta_j,v]\pa_xf\|_{L^{p}}\g\|_{\ell^{\infty}} \leq C\f(\|\pa_x v\|_{L^{\infty}}\|f\|_{B_{p, \infty}^{s}}+\|\pa_x f\|_{L^{\infty}}\|\pa_xv\|_{B_{p, \infty}^{s-1}}\g),
\end{align*}
where we denote the standard commutator $[\Delta_j,v]\pa_xf=\Delta_j(v\pa_xf)-v\Delta_j\pa_xf$.
\end{lemma}

\section{Proof of Theorem \ref{th1}}\label{sec3}
\subsection{Construction of Initial Data}\label{sec3.1}
Let $\widehat{\phi}\in \mathcal{C}^\infty_0(\mathbb{R})$ be an even, real-valued and non-negative function on $\R$ and satisfy
\begin{numcases}{\widehat{\phi}(\xi)=}
1,&if $|\xi|\leq \frac{1}{4}$,\nonumber\\
0,&if $|\xi|\geq \frac{1}{2}$.\nonumber
\end{numcases}
Obviously,  $\phi(0)>0$ and
for any $p\in[1,\infty]$, there exists two positive constants $c_1$ and $c_2$ such that
\begin{align*}
c_1\leq\|\phi\|_{L^p(\R)}\leq c_2.
\end{align*}
We
define the function $f_n(x)$ by
$$f_n(x)=\phi(x)\cos \f(\frac{17}{12}2^{n}x\g)\quad\text{with}\quad n\gg1.$$
Due to the fact $\varphi(\xi)\equiv 1$ for $\frac43\leq |\xi|\leq \frac32$,
we have
\begin{align}\label{cl}
{\Delta_j(f_n)=\mathcal{F}^{-1}\f(\varphi(2^{-j}\cdot)\widehat{f}_{n}\g)=}\begin{cases}
f_n, &\text{if}\; j=n,\\
0, &\text{if}\; j\neq n.
\end{cases}
\end{align}

\begin{lemma}\label{le5} Assume that $(s,p,r)$ satisfies \eqref{cond}.
Define the initial data $u_0(x)$ as
\bal\label{lyz-u0}
u_0(x):=\sum\limits^{\infty}_{n=3}n^{-2}2^{-ns} \phi(x)\cos \f(\frac{17}{12}2^{n}x\g).
\end{align}
Then there exists some sufficiently large $n\in \mathbb{Z}^+$  such that
\bbal
&\|u_0\|_{B^{s}_{p,r}}\leq C\quad\text{and}\quad \|u_0\pa_x\De_{n}u_0\|_{L^p}\geq cn^{-2}2^{n(1-s)},
\end{align*}
where $C$ and $c$ are some positive constants.
\end{lemma}

\begin{proof}\;
Using \eqref{cl} yields
\begin{align}\label{u5}
\De_{n}u_0(x)&=n^{-2}2^{-ns} \phi(x)\cos \f(\frac{17}{12}2^{n}x\g).
\end{align}
By the definition of ${B}_{p,r}^{s}$, we deduce that
\begin{align*}
 \|u_{0}\|_{{B}_{p,r}^{s}(\R)}&=\f\|2^{js}\|\Delta_{j}u_0\|_{L^p(\R)}\g\|_{\ell^r(j\geq1)}
\leq \f\|j^{-2}\g\|_{\ell^r(j\geq1)} \|\phi\|_{L^{p}(\R)}\leq C\|\phi\|_{L^{p}(\R)}.
\end{align*}
From \eqref{u5}, we have
\begin{align*}
n^2u_0\pa_x\De_{n}u_0&=2^{-ns} u_0(x)\phi'(x)\cos \f(\frac{17}{12}2^{n}x\g)-\frac{17}{12}2^{n}2^{-ns} u_0(x)\phi(x)\sin \f(\frac{17}{12}2^{n}x\g).
\end{align*}
Since $u_0(x)$ is a real-valued and continuous function on $\R$, then there exists some $\delta>0$ such that for any $x\in B_{\delta}(0)$
\begin{align}\label{yh}
&|u_0(x)|\geq \fr{1}{2}|u_0(0)|=\fr{1}{2}\phi(0)\sum\limits^{\infty}_{n=3}n^{-2}2^{-ns}=:c_0.
\end{align}
Thus we have from \eqref{yh}
\begin{align*}
n^2\|u_0\pa_x\De_{n}u_0\|_{L^p}
&\geq c_02^{n}2^{-ns} \f\|\phi(x)\sin \f(\frac{17}{12}2^{n}x\g)\g\|_{L^p(B_{\delta}(0))}-C2^{-ns}\f\| \phi'(x)\phi(x)\cos \f(\frac{17}{12}2^{n}x\g)\g\|_{L^p}\\
&\geq (c2^{n}-C)2^{-ns} .
\end{align*}
We choose $n$ large enough such that $C<\frac{c}{2}2^{n}$ and then finish the proof of Lemma \ref{le5}.
\end{proof}
\subsection{Error Estimates}\label{sec3.2}
\begin{proposition}\label{pro3.1}
Assume that $u_0$ satisfies \eqref{lyz-u0}. Under the assumptions of Theorem \ref{th1}, we have
\begin{align*}
&\|\mathbf{S}_{t}\left(u_{0}\right)-u_0\|_{B^{s-1}_{p,r}}\lesssim t.
\end{align*}
Furthermore, there holds
\begin{itemize}
  \item for $s>\max\f\{1+\frac{1}{p},\frac{3}{2}\g\}$ and $(p, r)\in [1, \infty]\times [1, \infty),$ we have
  \begin{align*}
\|\mathbf{w}\|_{B^{s-2}_{p,r}}\lesssim t^{2},
\end{align*}
here and in what follows we denote
$\mathbf{w}:=\mathbf{S}_{t}(u_0)-u_0-t\mathbf{\widetilde{u}}_0$ with $\mathbf{\widetilde{u}}_0:=\mathbf{P}(u_0)-u_0\pa_x u_0.$
  \item for $s\in\big(1+\frac{1}{p},\frac{3}{2}\big]$ and $(p, r)\in [1, \infty]\times [1, \infty)$ or $s=1+\frac{1}{p}$ and $(p, r)\in [1, \infty)\times \{1\}$, we have
  \begin{align*}
\|\mathbf{w}\|_{B^{0}_{p,r}}\lesssim t^{s}.
\end{align*}
\end{itemize}
\end{proposition}
\begin{proof}\; For simplicity, we denote $u(t):=\mathbf{S}_t(u_0)$ here and in what follows. Notice that $(s,p,r)$ satisfies \eqref{cond} and by the local well-posedness result \cite{d1,Li-Yin1,Ye}, there exists a positive time $T$ such that $u(t)\in \mathcal{C}([0,T];B_{p,r}^s)$. Furthermore, it holds that
$
\|u(t)\|_{L^\infty_TB^s_{p,r}}\leq C\|u_0\|_{B^s_{p,r}}\leq C.
$
For any $u_0\in B^s_{\infty,r}$, we can deduce that $u(t)\in \mathcal{C}([0,T];B_{\infty,r}^s)$ and
$
\|u(t)\|_{L^\infty_TB^s_{\infty,r}}\leq C\|u_0\|_{B^s_{\infty,r}}\leq C.
$
Using the Newton-Leibniz formula and Remark \ref{re5} with $v=0$, we obtain from \eqref{CH} that
\bal\label{s}
\|u(t)-u_0\|_{B^{s-1}_{p,r}}
&\leq \int^t_0\|\pa_\tau u\|_{B^{s-1}_{p,r}} \dd\tau
\leq \int^t_0\|\mathbf{P}(u)\|_{B^{s-1}_{p,r}}+\|u\pa_xu\|_{B^{s-1}_{p,r}} \dd\tau
\les t\|u\|^{2}_{L_t^\infty B^{s}_{p,r}}
\les t\|u_0\|^{2}_{B^{s}_{p,r}}\les t.
\end{align}
{\bf Case 1:} $s>\max\f\{1+\frac{1}{p},\frac{3}{2}\g\}$ and $(p, r)\in [1, \infty]\times [1, \infty)$.
By the Newton-Leibniz formula and Lemmas \ref{lm2}-\ref{lm3}, we obtain from \eqref{s} that
\begin{align*}
\|\mathbf{w}\|_{B^{s-2}_{p,r}}
\leq &~ \int^t_0\|\partial_\tau u-\mathbf{\widetilde{u}}_0\|_{B^{s-2}_{p,r}} \dd\tau
\les \int^t_0\|\mathbf{P}(u)-\mathbf{P}(u_0)\|_{B^{s-2}_{p,r}} +\|u\partial_xu-u_0\partial_xu_0\|_{B^{s-2}_{p,r}}\dd\tau
\nonumber\\
\les &~\int^t_0\|\pa_x(u-u_0)\pa_x(u+u_0)\|_{B^{s-2}_{p,r}}+\|(u-u_0)(u+u_0)\|_{B^{s-1}_{p,r}}\dd\tau\nonumber
\\
\les &~ \int^t_0\|\pa_x(u-u_0)\|_{B^{s-2}_{p,r}}\|\pa_x(u+u_0)\|_{B^{s-1}_{p,r}} +\|u-u_0\|_{B^{s-1}_{p,r}}\|u+u_0\|_{B^{s-1}_{p,r}}\dd\tau
\les  t^2.
\end{align*}
{\bf Case 2:} $1+\frac{1}{p}<s\leq\frac{3}{2}$ and $(p, r)\in [1, \infty]\times [1, \infty)$ or $s=1+\frac{1}{p}$ and $(p, r)\in [1, \infty)\times \{1\}$.
Using the embedding $B^{s-2}_{p,r}\hookrightarrow B^{-1}_{p,r}$, we have
\begin{align}\label{p2}
\|\mathbf{w}\|_{B^{0}_{p,r}}
\leq &~ \int^t_0\|\partial_\tau u-\mathbf{\widetilde{u}}_0\|_{B^{0}_{p,r}} \dd\tau
\les \int^t_0\|\mathbf{P}(u)-\mathbf{P}(u_0)\|_{B^{0}_{p,r}} +\|u\partial_xu-u_0\partial_xu_0\|_{B^{0}_{p,r}}\dd\tau
\nonumber\\
\les &~ \int^t_0\|\pa_x(u-u_0)\pa_x(u+u_0)\|_{B^{-1}_{p,r}} +\|(u-u_0)(u+u_0)\|_{B^{1}_{p,r}}\dd\tau\nonumber
\\
\les &~ \int^t_0\underbrace{\|\pa_x(u(\tau)-u_0)\pa_x(u(\tau)+u_0)\|_{B^{s-2}_{p,r}}}_{=:\ I_1(\tau)} +\underbrace{\|u(\tau)-u_0\|_{B^{1}_{p,r}}\|u(\tau)+u_0\|_{B^{1}_{p,r}}}_{=:\ I_2(\tau)}\dd\tau.
\end{align}
To estimate the term $I_1$, by Bony's decomposition, one has
\bbal
\pa_x(u(t)-u_0)\pa_x(u(t)+u_0)&=T_{\pa_x(u(t)-u_0)}\pa_x(u(t)+u_0)+T_{\pa_x(u(t)+u_0)}\pa_x(u(t)-u_0)+R\big(\pa_x(u(t)+u_0),\pa_x(u(t)-u_0)\big).
\end{align*}
Using Lemma \ref{leem1} yields that
\bbal
\|T_{\pa_x(u(t)-u_0)}\pa_x(u(t)+u_0)\|_{B^{s-2}_{p,r}}&\les \|\pa_x(u(t)-u_0)\|_{B^{-1}_{\infty,\infty}}\|\pa_x(u(t)+u_0)\|_{B^{s-1}_{p,r}}
\\&\les \|u(t)-u_0\|_{B^{1}_{p,r}}\|u(t)+u_0\|_{B^{s}_{p,r}}\les t^{s-1},\\
\|T_{\pa_x(u(t)+u_0)}\pa_x(u(t)-u_0)\|_{B^{s-2}_{p,r}}&\les \|\pa_x(u(t)+u_0)\|_{L^\infty}\|\pa_x(u(t)-u_0)\|_{B^{s-2}_{p,r}}\\&\les \|u(t)+u_0\|_{B^{s-1}_{p,r}}\|u(t)-u_0\|_{B^{s-1}_{p,r}}\les t,\\
\|R\big(\pa_x(u(t)+u_0),\pa_x(u(t)-u_0)\big)\|_{B^{s-1}_{p,r}}&\les \|\pa_x(u(t)+u_0)\|_{B^{s-1}_{\infty,\infty}}\|\pa_x(u(t)-u_0)\|_{B^{0}_{p,r}}
\\&\les \|u(t)+u_0\|_{B^{s}_{\infty,r}}\|u(t)-u_0\|_{B^{1}_{p,r}}\les t^{s-1},
\end{align*}
where we have used the fact
$\|u(t)\|_{B^{s}_{\infty,r}}\les \|u_0\|_{B^{s}_{\infty,r}}\les 1$
and  the interpolation argument
$$\|u(t)-u_0\|_{B^{1}_{p,r}}\les \|u(t)-u_0\|^{s-1}_{B^{s-1}_{p,r}}\|u(t)-u_0\|^{2-s}_{B^{s}_{p,r}}\les t^{s-1}.$$
Combining the above, we obtain
\bal\label{j1}
I_1(t)\les t^{s-1}.
\end{align}
To estimate the term $I_2$, by the interpolation argument, one has
\bal\label{j2}
I_2(t)\les \|u(t)-u_0\|_{B^{1}_{p,r}}\les t^{s-1}.
\end{align}
Inserting \eqref{j1} and \eqref{j2} into \eqref{p2} yields the desired result.
Thus, we finish the proof of Proposition \ref{pro3.1}.
\end{proof}

Now we present the proof of Theorem \ref{th1}.
\begin{proof}

{\bf Case 1:} $s>\max\f\{1+\frac{1}{p},\frac{3}{2}\g\}$ and $(p, r)\in [1, \infty]\times [1, \infty)$.
Notice that $\mathbf{S}_{t}(u_0)-u_0=t\mathbf{\widetilde{u}}_0+\mathbf{w}$ and $\mathbf{\widetilde{u}}_0=\mathbf{P}(u_0)-u_0\pa_x u_0.$ By the triangle inequality and Proposition \ref{pro3.1}, we deduce that
\bal\label{M}
\|\mathbf{S}_{t}(u_0)-u_0\|_{B^{s}_{p,r}}
&\geq2^{{ns}}
\big\|\De_{n}\big(\mathbf{S}_{t}(u_0)-u_0\big)\big\|_{L^p}
=2^{{ns}}\big\|\De_{n}\big(t\mathbf{\widetilde{u}}_0+\mathbf{w}\big)\big\|_{L^p}
\nonumber\\&\geq t2^{{ns}}\|\De_{n}\mathbf{\widetilde{u}}_0\|_{L^p}
-2^{{2n}}2^{{n(s-2)}}
\big\|\De_{n}\mathbf{w}\big\|_{L^p}\nonumber\\
&\geq t2^{{n}s}\|\De_{n}\big(u_0\pa_xu_0\big)\|_{L^p}-
t2^{{n}s}\|\De_{n}\big(\mathbf{P}(u_0)\big)\|_{L^p}-C2^{2{n}}\|\mathbf{w}\|_{B^{s-2}_{p,\infty}}
\nonumber\\&\geq t2^{{n}s}\|u_0\pa_x\De_{n}u_0\|_{L^p}-t2^{{n}s}\|[\De_{n},u_0]\pa_xu_0\|_{L^p}-
Ct\|\mathbf{P}(u_0)\|_{B^{s}_{p,\infty}}-C2^{2{n}}t^2\nonumber\\&\geq t2^{{n}s}\|u_0\pa_x\De_{n}u_0\|_{L^p}-Ct\big\|2^{{n}s}\|[\De_{n},u_0]\pa_xu_0\|_{L^p}\big\|_{\ell^\infty}-t\|\mathbf{P}(u_0)\|_{B^{s}_{p,\infty}}-C2^{2{n}}t^2.
\end{align}
By Lemmas \ref{lm2}-\ref{lm5}, one has
\bbal
&\|\mathbf{P}(u_0)\|_{B^{s}_{p,\infty}}\les \f\|u_0^{2}+\fr12(\pa_xu_0)^2\g\|_{B^{s-1}_{p,r}}\les 1,\\
&\big\|2^{{n}s}\|[\De_{n},u_0]\pa_xu_0\|_{L^p}\big\|_{\ell^\infty}\les \|\pa_xu_0\|_{L^\infty}\|u_0\|_{B^{s}_{p,\infty}}+
\|\pa_xu_0\|_{L^\infty}\|\pa_xu_0\|_{B^{s-1}_{p,\infty}}\les 1.
\end{align*}
Gathering all the above estimates and Lemma \ref{le5} together with \eqref{M}, we obtain
\bbal
\|\mathbf{S}_{t}(u_0)-u_0\|_{B^{s}_{p,r}}\geq ctn^{-2}2^{n}-Ct-C2^{2{n}}t^2,
\end{align*}
which implies
\bbal
t^{-\alpha}\|\mathbf{S}_{t}(u_0)-u_0\|_{B^{s}_{p,r}}\geq ct^{1-\alpha}n^{-2}2^{n}-Ct^{1-\alpha}-C2^{2{n}}t^{2-\alpha}.
\end{align*}
Thus, picking $t^{1-\alpha}_n=n^32^{-n}$ with large $n$, we have
\bbal
t^{-\alpha}_n\|\mathbf{S}_{t_n}(u_0)-u_0\|_{B^{s}_{p,r}}&\geq cn-Cn^32^{-n}-Cn^{6}t_n^{\alpha}
\geq \tilde{c}n.
\end{align*}
{\bf Case 2:} $1+\frac{1}{p}<s\leq\frac{3}{2}$ and $(p, r)\in [1, \infty]\times [1, \infty)$ or $s=1+\frac{1}{p}$ and $(p, r)\in [1, \infty)\times \{1\}$.
By Proposition \ref{pro3.1}, we deduce that
\bal\label{M2}
\|\mathbf{S}_{t}(u_0)-u_0\|_{B^{s}_{p,r}}
&\geq2^{{ns}}
\big\|\De_{n}\big(\mathbf{S}_{t}(u_0)-u_0\big)\big\|_{L^p}=2^{{ns}}\big\|\De_{n}\big(t\mathbf{\widetilde{u}}_0+\mathbf{w}\big)\big\|_{L^p}
\nonumber\\&\geq t2^{{ns}}\|\De_{n}\mathbf{\widetilde{u}}_0\|_{L^p}
-2^{{ns}}
\|\mathbf{w}\|_{B^0_{p,r}}\geq t2^{{n}s}\|\De_{n}\big(u_0\pa_xu_0\big)\|_{L^p}-
t2^{{n}s}\|\De_{n}\big(\mathbf{P}(u_0)\big)\|_{L^p}-C(2^{{n}}t)^{s}
\nonumber\\&\geq t2^{{n}s}\|u_0\pa_x\De_{n}u_0\|_{L^p}-t2^{{n}s}\|[\De_{n},u_0]\pa_xu_0\|_{L^p}-
t\|\mathbf{P}(u_0)\|_{B^{s}_{p,r}}-C(2^{{n}}t)^{s}\nonumber\\&\geq t2^{{n}s}\|u_0\pa_x\De_{n}u_0\|_{L^p}-Ct\big\|2^{{n}s}\|[\De_{n},u_0]\pa_xu_0\|_{L^p}\big\|_{\ell^\infty}-t\|\mathbf{P}(u_0)\|_{B^{s}_{p,\infty}}-C(2^{{n}}t)^{s},
\end{align}
where we have used
\bbal
&\|\mathbf{P}(u_0)\|_{B^{s}_{p,r}}\les \f\|u_0^{2}+\fr12(\pa_xu_0)^2\g\|_{B^{s-1}_{p,r}}\les 1,\\
&\big\|2^{{n}s}\|[\De_{n},u_0]\pa_xu_0\|_{L^p}\big\|_{\ell^\infty}\les \|\pa_xu_0\|_{L^\infty}\|u_0\|_{B^{s}_{p,\infty}}+
\|\pa_xu_0\|_{L^\infty}\|\pa_xu_0\|_{B^{s-1}_{p,\infty}}\les 1.
\end{align*}
Gathering Lemma \ref{le5} together with \eqref{M2}, we obtain
\bbal
\|\mathbf{S}_{t}(u_0)-u_0\|_{B^{s}_{p,r}}\geq ctn^{-2}2^{{n}}-Ct-C(2^{{n}}t)^{s},
\end{align*}
which implies
\bbal
t^{-\alpha}\|\mathbf{S}_{t}(u_0)-u_0\|_{B^{s}_{p,r}}\geq ct^{1-\alpha}n^{-2}2^{{n}}-Ct^{1-\alpha}-C(2^{{n}}t)^{s}t^{-\alpha}.
\end{align*}
Thus, picking $t^{1-\alpha}_n=n^32^{-n}$ with large $n$, we have
\bbal
t^{-\alpha}_n\|\mathbf{S}_{t_n}(u_0)-u_0\|_{B^{s}_{p,r}}&\geq cn-Cn^32^{-n}-Cn^{3s}t_n^{\alpha(s-1)}
\geq \tilde{c}n.
\end{align*}
This completes the proof of Theorem \ref{th1}.\end{proof}

\section*{Acknowledgments}
The authors would like to thank the anonymous referees for valuable comments and suggestions which greatly improved the presentation of this paper. J. Li is supported by Jiangxi Provincial Natural Science Foundation (20252BAB22004). W. Zhu is supported by the National Natural Science Foundation of China (12201118).

\end{document}